\documentclass[12pt]{article}
\usepackage{amssymb,amsfonts,amsmath, psfrag,eepic,colordvi}
\topskip  
\numberwithin{equation}{section}
\newtheorem{theorem}{Theorem}[section]

\begin{document}
\parskip 6pt

\pagenumbering{arabic}
\def\sof{\hfill\rule{2mm}{2mm}}
\def\ls{\leq}
\def\gs{\geq}
\def\SS{\mathcal S}
\def\qq{{\bold q}}
\def\MM{\mathcal M}
\def\TT{\mathcal T}
\def\EE{\mathcal E}
\def\lsp{\mbox{lsp}}
\def\rsp{\mbox{rsp}}
\def\pf{\noindent {\it Proof.} }
\def\mp{\mbox{pyramid}}
\def\mb{\mbox{block}}
\def\mc{\mbox{cross}}
\def\qed{\hfill \rule{4pt}{7pt}}
\def\block{\hfill \rule{5pt}{5pt}}

\begin{center}
{\Large\bf    On   partitions avoiding right crossings
  } \vskip 6mm
\end{center}

\begin{center}
{\small   Sherry H. F. Yan,    Yuexiao Xu \\[2mm]
 Department of Mathematics, Zhejiang Normal University, Jinhua
321004, P.R. China
\\[2mm]
 huifangyan@hotmail.com
  \\[0pt]
}
\end{center}

\noindent {\bf Abstract.} Recently, Chen et al. derived the
generating function for partitions avoiding right nestings and posed
the problem of finding the generating function for partitions
avoiding right crossings.  In this paper, we derive the generating
function for partitions avoiding right crossings via an intermediate
structure of partial matchings avoiding $2$-right crossings and
right nestings. We show that there is a bijection between partial
matchings avoiding 2-right crossing and right nestings and
partitions avoiding right crossings.




\section{Introduction}

This paper is concerned with the enumeration of   set partitions
that avoid right crossings. Recall that a partition $\pi$ of
$[n]=\{1,2,\ldots, n\}$ can be represented as a diagram with
vertices drawn on a horizontal line in increasing order. For a block
$B$ of $\pi$, we write the elements of $B$ in increasing order.
Suppose that $B=\{i_1, i_2, \ldots, i_k\}$. Then we draw an arc from
$i_1$ to $i_2$, an arc from $i_2$ to $i_3$, and so on. Such a
diagram is called the {\em linear representation} of $\pi$.  If
$(i,j)$ is an arc in the diagram of $\pi$ with $i<j$, we call vertex
$i$ an {\em opener}, and call vertex $j$ a {\em transient} if there
is an arc $(j,k)$ for some integer $k$ such that $i<j<k$; otherwise,
we call vertex $j$ a {\em closer}. A partial matching is a partition
for which each block contains at most two elements. A partial
matching without singletons is called a perfect matching.

A {\em nesting} of a partition is a pair of arcs $(i_1, j_1)$ and
$(i_2, j_2)$ with $i_1<i_2<j_2<j_1$. We call such a nesting a {\em
left nesting} if $i_2=i_1+1$. Similarly, we call it a {\em right
nesting} if $j_1=j_2+1$. A {\em crossing} of a partition is a pair
of arcs $(i_1, j_1)$ and $(i_2, j_2)$ with $i_1<i_2<j_1<j_2$ and we
can define {\em left crossing } and {\em right crossing }
analogously to how it was defined for nesting arcs.  A {\em neighbor
alignment} of a partition is a pair of arcs $(i_1,j_1)$ and $(i_2,
j_2)$ such that $i_2=j_1+1$.

 Zagier
\cite{zag} derived  the generating function for perfect  matchings
avoiding left nestings and right nestings, which was  introduced by
Stoimenow \cite{Stoi}. Recently,  Bousquet-M\'elou et al.
\cite{melon} constructed bijections between perfect matchings
avoiding left nestings and right nestings and three other classes of
combinatorial objects, that is,  unlabeled $(2+2)$-free posets,
permutations avoiding a certain pattern and ascent sequences. Dukes
and Parviainen \cite{Duck} presented a bijection between ascent
sequences and upper triangular matrices whose non-negative entries
are such that all rows and columns contain at least one non-zero
entry.  Claesson and Linusson \cite{cla} established a bijection
between perfect matchings avoiding left nestings and permutations
and showed that perfect matchings avoiding left crossings are
equinumerous with perfect matchings avoiding left nestings.
 They also introduced the notion of {\em k-left
nesting} and conjectured that the number of
 perfect matchings of $[2n]$ without 2-left nestings
 is equal to  the number of perfect matchings of $[2n]$ without left nestings and right
 nestings.  This conjecture was confirmed by
Levande \cite{levande}.  Recall that two arcs $(i_1, j_1)$ and
$(i_2, j_2)$ form
 a {\em $k$-left nesting} if $i_1<i_2<j_2<j_1$ and $i_2-i_1\leq k$.  Similarly, we call it a {\em $k$-right
nesting} if $i_1<i_2<j_2<j_1$ and $j_1-j_2\leq k$. Note that  1-left
 nesting is exactly a left nesting and 1-right
 nesting is exactly a right nesting. Similarly,  we
can define {\em k-left crossing } and {\em k-right crossing }
analogously to how it was defined for nesting arcs. The left
nestings, left crossings, right crossings, right crossings, neighbor
alignments are called neighbor patterns.

Recently, Chen et al. \cite{chen}  derived the generating functions
for partial matchings avoiding neighbor alignments and for partial
matchings avoiding neighbor alignments and left nestings. They
obtained the generating function for partitions avoiding right
nestings by   presenting  a bijection between partial matchings
avoiding three neighbor patterns ( left nestings, right nestings and
neighbor alignments) and partitions avoiding right nestings. In
general, the number of partitions of $[n]$ avoiding right crossings
is not equal to the number of partitions of $[n]$ avoiding right
nestings.  In this paper, we derive the generating functions  for
partitions avoiding right crossings  via an intermediate structure
of partial matchings avoiding $2$-right crossings and right
nestings.

 Denote by $\mathcal{M}(n,k)$ the set of partial matchings
of $[n]$ with $k$
 arcs.
 The set of partial matchings in
$\mathcal{M}(n,k)$ with no $2$-right crossings   and right nestings
is denoted by $\mathcal{P}(n,k)$.    Let $\mathcal{T}(n,k)$ be the
set of partitions of $[n]$ with $k$ arcs. Denote by $\mathcal{CT}(
n,k)$ the
 set of partition in $\mathcal{T}(n,k)$
 with no right crossings.

 Denote by $P(n,k)$   and  $CT(n,k)$ the
cardinalities of the sets $\mathcal{P}(n,k)$  and
$\mathcal{CT}(n,k)$, respectively.

We obtain the  generating function formula for the numbers
$P(n+k-1,k)$ by  establishing a bijection between  partial matchings
avoiding 2-right crossings and right nestings and a certain class of
integer sequence. Moreover, we show that there is a correspondence
between   $\mathcal{P}(n+k-1,k)$ and   $\mathcal{CT}(n,k)$.

\begin{theorem}\label{thl} We have
\begin{equation}\label{L(n,k)}
\sum_{n\geq 1}\sum_{k=0}^{n-1}P(n+k-1,k)x^ny^k=\sum_{n\geq 1}
{{x^n(1+xy)^{n\choose 2}}\over  \prod_{k=0}^{n-1}(1-(1+xy)^kxy) }.
\end{equation}
\end{theorem}

\begin{theorem}\label{tht} We have
\begin{equation}\label{T(n,k)}
\sum_{n\geq 1}\sum_{k=0}^{n-1}CT(n,k)x^ny^k=\sum_{n\geq 1}
{{x^n(1+xy)^{n\choose 2}}\over  \prod_{k=0}^{n-1}(1-(1+xy)^kxy) }.
\end{equation}
\end{theorem}

\section{ Partial matchings avoiding 2-right crossings and right nestings }

In this section, we construct a bijection between  partial matchings
avoiding 2-right crossings and right nestings and a certain class of
integer sequences. As a consequence, we obtain the bivariate
generating function for the number of partial matchings of $[n+k-1]$
with $k$ arcs and containing  no 2-right crossing and right
nestings.

 Let $x=x_1x_2\ldots x_n$ be an integer sequence. Denote by
$max(x)$ the maximum element of $x$. For $1<i\leq n$, an element
$x_i$ of $x$ is said to be a {\em  left-to-right maximum} if $x_i>
max(x_1x_2\ldots x_{i-1})$.  For $n\geq 1$, let $\mathcal{S}(n,k)$
be the set of integer sequences $x=x_0x_1\ldots x_{n-1}$ with $n-k$
left-to-right maxima satisfying that
\begin{itemize}
\item $x_0= 0$;
\item  for all $1\leq i\leq n-1$,  $0\leq x_i\leq max(x_0x_1\ldots x_{i-1})+1$;
\item  for all $0\leq i\leq n-1$, if $0\leq  x_i< max(x_0x_1\ldots x_{i-1})$, then $x_i<x_{i-1}$.
\end{itemize}
Denote by $\mathcal{S}(n)$ the set of such integer sequences  with
$n$ left-to-right maxima.   Denote by  $S(n,k)$  the cardinality of
the set $\mathcal{S}(n,k)$.  Let $f(x,y)$ be the generating function
for the numbers $S(n,k)$. We derive the following generating
function formula of $f(x,y)$.

\begin{theorem}\label{ths}
\begin{equation}\label{S(n,k)}
f(x,y)=\sum_{n\geq 1}\sum_{k=0}^{n-1}S(n,k)x^ny^k=\sum_{n\geq 1}
{{x^n(1+xy)^{n\choose 2}}\over  \prod_{k=0}^{n-1}(1-(1+xy)^kxy) }.
\end{equation}
\end{theorem}
\pf Let $x$ be a sequence in $\mathcal{S}(n)$ with $n$ left-to-right
maxima.
  It can be uniquely decomposed as
$0w_01w_12w_2\ldots  (n-1)w_{n-1}$. For all $0\leq i\leq n-1$, each
$w_i$ is a (possibly empty) integer sequence whose elements are
nonnegative integers less than or equal to $i$. For all $1\leq i\leq
n-1$, suppose that there are $i_l$ occurrences of $i$'s in $w_i$,
then $w_i$ reads
 as $m_0im_1im_2\ldots i m_{i_l}$,  where for all $0\leq j\leq
i_l$, the sequence  $m_j$ is a (possibly empty) decreasing sequence
whose elements are nonnegative integers less than $i$. Since  $w_0$
is a (possibly empty) sequence of $0$'s. Hence, $w_0$ contributes
$1\over {1-xy}$ to the generating function $f(x,y)$, while for all
$1\leq i\leq n-1$, $w_i$ contributes $(1+xy)^i\over {1-(1+xy)^ixy}$
to the generating function $f(x,y)$. Furthermore, each left-to-right
maximum  contributes $x$ to the generating function $f(x,y)$.
Summing over all $n\geq 1$, we derive  the generating function
formula (\ref{S(n,k)}). This completes the proof.\qed

Now we proceed to construct a bijection between the set
$\mathcal{P}(n+k-1,k)$ and the set $\mathcal{S}(n,k)$.
\begin{theorem}\label{bijection}
There exists a bijection  between the set $\mathcal{P}(n+k-1,k)$ and
the set $\mathcal{S}(n,k)$.
\end{theorem}
 \pf Let $M$ be a partial matching in $\mathcal{P}(n+k-1,k)$ of
 $[n+k-1]$ with $k$ arcs and containing no 2-right crossings and
 right nestings. We wish to generate   a sequence   $\alpha(M)=x_0x_1x_2\ldots x_{n-1}$ from $M$ recursively. First,
  we remove the labels of the closers and relabel the vertices of $M$ in the natural order.
  Denote by $\overline{M}$ the obtained diagram.
  For all $1\leq i\leq n-1$, let  $\mathcal{O}(i)$
be the set of openers of arcs whose closers are  left  to vertex $i$
of $\overline{M}$.  Denote by $\mathcal{O}(n)$  the set of openers
of all arcs of $\overline{M}$.
   Set $x_0=0$ and  assume that we have
 obtained $x_{i-1}$. Now we proceed to generate $x_i$ in the
 following manner.
 \begin{itemize}
 \item[Case 1.]
  If there is no closer immediately after vertex $i$, then let $x_i=max(x_0x_1\ldots
 x_{i-1})+1$.
\item[Case 2.] If  there is an arc whose closer is immediately after vertex
$i$ and opener is labelled with  $j $ and $j$ is the $m$-th minimum
element  of the set $[i]\setminus\mathcal{O}(i)$, then let
$x_i=m-1$.
\end{itemize}

We claim that  for all $0\leq i\leq n-1$, $max(x_0x_1x_2\ldots x_i)=
i-|\mathcal{O}(i+1)|$. We prove this claim by induction on $i$.
Obviously, the claim holds for $i=0$.  So let $i\geq  1$ and assume
that the claim holds for $i-1$. In Case 1, we have
$max(x_0x_1x_2\ldots x_i)=max(x_0x_1\ldots
x_{i-1})+1=i-1-|\mathcal{O}(i)|+1$. Since there exist no closers
immediately after vertex $i$ in $\overline{M}$, we are led to
$i-|\mathcal{O}(i+1) |=i-1-|\mathcal{O}(i)| +1=max(x_0x_1x_2\ldots
x_i)$. In Case 2,  we have $0\leq x_i\leq
i-|\mathcal{O}(i)|-1=max(x_0x_1\ldots x_{i-1})$, which yields that
$max(x_0x_1\ldots x_{i})=max(x_0x_1\ldots
x_{i-1})=i-1-|\mathcal{O}(i)|$. Since there is exactly one closer
immediately after vertex $i$, we have
$i-|\mathcal{O}(i+1)|=i-1-|\mathcal{O}(i)| =max(x_0x_1x_2\ldots
x_i)$. Therefore, the claim holds for all $0\leq i\leq n-1$.   This
yields that
 $0\leq x_i\leq max(x_0x_1x_2\ldots x_{i-1})+1$ for all $1\leq i\leq
 n-1$ and
$max(x_0x_1\ldots x_{n-1})=n-1-|\mathcal{O}(n)|=n-1-k$, which
implies that there are $n-k$ left-to-right maxima in $\alpha(M)$.

In order to prove that the obtained sequence $\alpha(M)\in
\mathcal{S}(n,k)$, it remains to show that in  the sequence
$\alpha(M)$  if $0\leq x_i< max(x_1x_2\ldots x_{i-1})$, then
$x_i<x_{i-1}$. According to the construction of the map $\alpha$,
there exists  an arc whose closer is immediately after vertex $i$.
Suppose that the opener of this arc is labelled with $p$ in
$\overline{M}$. Since $x_i< max(x_1x_2\ldots
x_{i-1})=i-1-|\mathcal{O}(i)|$, we are led to $p<i$.  Obviously,
when $x_{i-1}=max(x_1x_2\ldots x_{i-1})$, we have $x_i<x_{i-1}$. Now
suppose that $x_{i-1}<max(x_1x_2\ldots x_{i-1})$. According to the
construction of the map $\alpha$, there exists  an arc whose closer
is immediately after vertex $i-1$. Suppose that the opener of this
arc is labelled with $q$ in $\overline{M}$. Since there are no
2-right crossings in $M$, we have $p<q$. It is easy to check that
$\mathcal{O}(i)=\mathcal{O}(i-1)\cup \{q\}$. From the construction
of the map $\alpha$, we have  $x_i<x_{i-1}$. Hence, the obtained
sequence $\alpha(M)\in \mathcal{S}(n,k)$.

Conversely, given a sequence  $x=x_0x_1x_2\ldots x_{n-1}\in
\mathcal{S}(n, k)$, we wish to construct a partial matching of
$[n+k-1]$ with $k$ arcs. First, we arrange $n-1$ vertices on a
horizontal line and label them in the natural order. For
$i=1,2,\ldots, n-1$, at step $i$, we insert at most one closer after
vertex $i$ described as follows:
\begin{itemize}
\item[(i)] if $x_i=max(x_0x_1x_2\ldots x_{i-1})+1$,  then we do nothing for vertex
$i$;
 \item[(ii)] if $0\leq x_i\leq  max(x_0x_1x_2\ldots x_{i-1})$, then we insert  an arc whose opener is
the $(x_i+1)$-th vacant vertex (from left to right) of the set $[i]$
and closer is immediately after vertex $i$, where a vacant vertex is
a vertex which is not joined to any arc.
\end{itemize}
  Finally, we relabel the
vertices in the natural order to get a partial matching.

By induction on $i$ $(0\leq i\leq n-1)$, it is easy to see that
after step $i$, the number of  vacant vertices  of the set $[i]$ is
equal to $max(x_0x_1x_2\ldots x_{i})$. So (ii) is valid and the
construction of the inverse map is well defined. Since  there are
exactly $n-k$ left-to-right maxima in $x$,   the obtained partial
matching is a matching of $[n+k-1]$ with $k$ arcs.  Note that we add
at most one closer immediately after each vertex $i$ for all $1\leq
i\leq n-1$. Thus there is no two consecutive closers in the
resulting partial matchings, which implies that there are no two
right crossings and right nestings. Since when $0\leq
x_i<max(x_1x_2\ldots x_{i-1})$, we have $x_i<x_{i-1}$, there are no
2-right crossings. Hence, the resulting partial matching is in the
set $\mathcal{P}(n+k-1, k)$. This implies that the above map
$\alpha$ is a bijection.      \qed

\begin{figure}[h,t]
\begin{center}
\begin{picture}(60,10)
\setlength{\unitlength}{6mm} \linethickness{0.4pt}
 \put(-0.2,-0.6){\small$1$}\put(0.8,-0.6){\small$2$}
\put(2.8,-0.6){\small$3$}\put(3.8,-0.6){\small$4$}
 \put(5.8,-0.6){\small$5$} \put(7.8,-0.6){\small$6$}
\put(9.8,-0.6){\small$7$}
\qbezier(0,0)(2.5,2.5)(5,0)\qbezier(1,0)(1.5,1.5)(2,0)\qbezier(3,0)(7,3)(11,0)\qbezier(4,0)(6.5,2.5)(9,0)
\qbezier(6,0)(6.5,1.5)(7,0)
\put(8,0){\circle*{0.2}}\put(10,0){\circle*{0.2}}\put(0,0){\circle*{0.2}}\put(1,0){\circle*{0.2}}\put(2,0){\circle*{0.2}}
\put(3,0){\circle*{0.2}}\put(4,0){\circle*{0.2}}\put(5,0){\circle*{0.2}}\put(6,0){\circle*{0.2}}\put(7,0){\circle*{0.2}}
\put(9,0){\circle*{0.2}}\put(11,0){\circle*{0.2}}

\end{picture}
\vspace{-15pt}
\end{center}
\caption{The diagram $\overline{M}$.} \label{fig1}
\end{figure}
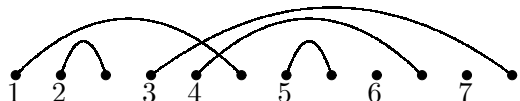

For example, let $M=\{\{1,6\}, \{2,3\}, \{4,12\}, \{5,10\}, \{7,8\},
\{9\},  \{11\} \}\in \mathcal{P}(12,5)$. We obtain a diagram
$\overline{M}$ from   $M$ by removing the labels of the closers and
relabeling the vertex of $M$ in the natural order, see Figure
\ref{fig1}. We have $\mathcal{O}(1)=\emptyset$,
$\mathcal{O}(2)=\emptyset$, $\mathcal{O}(3)=\{2\}$,
$\mathcal{O}(4)=\{2\}, \mathcal{O}(5)=\{1,2\},
\mathcal{O}(6)=\{1,2,5\}$, $\mathcal{O}(7)=\{1,2,4,5\}$ and
$\mathcal{O}(8)=\{1,2,3,4,5\}$. So the corresponding sequence  is
$\alpha(M)= 01120210$.

\begin{figure}[h,t]
\begin{center}
\begin{picture}(60,150)
\setlength{\unitlength}{6mm} \linethickness{0.4pt}
 \put(-0.2,0.4){\small$1$}\put(0.8,0.4){\small$2$}
\put(1.8,0.4){\small$3$}\put(2.8,0.4){\small$4$}
 \put(3.8,0.4){\small$5$} \put(4.8,0.4){\small$6$}
\put(5.8,0.4){\small$7$}\put(6.8,0.4){\small$8$}\put(7.8,0.4){\small$9$}\put(8.8,0.4){\small$10$}\put(9.8,0.4){\small$11$}
\put(10.8,0.4){\small$12$}
\qbezier(0,1)(2.5,3.5)(5,1)\qbezier(1,1)(1.5,2.5)(2,1)\qbezier(3,1)(7,4)(11,1)\qbezier(4,1)(6.5,3.5)(9,1)
\qbezier(6,1)(6.5,2.5)(7,1)
\put(8,1){\circle*{0.2}}\put(10,1){\circle*{0.2}}\put(0,1){\circle*{0.2}}\put(1,1){\circle*{0.2}}\put(2,1){\circle*{0.2}}
\put(3,1){\circle*{0.2}}\put(4,1){\circle*{0.2}}\put(5,1){\circle*{0.2}}\put(6,1){\circle*{0.2}}\put(7,1){\circle*{0.2}}
\put(9,1){\circle*{0.2}}\put(11,1){\circle*{0.2}}
\put(6,3.5){$\downarrow$} \put(-3, 5){\small $x_7=0$}
 \put(-0.2,4.4){\small$1$}\put(0.8,4.4){\small$2$}
\put(2.8,4.4){\small$3$}\put(3.8,4.4){\small$4$}
 \put(5.8,4.4){\small$5$} \put(7.8,4.4){\small$6$}
\put(9.8,4.4){\small$7$}
\qbezier(0,5)(2.5,7.5)(5,5)\qbezier(1,5)(1.5,6.5)(2,5)\qbezier(3,5)(7,8)(11,5)\qbezier(4,5)(6.5,7.5)(9,5)
\qbezier(6,5)(6.5,6.5)(7,5)
\put(8,5){\circle*{0.2}}\put(10,5){\circle*{0.2}}\put(0,5){\circle*{0.2}}\put(1,5){\circle*{0.2}}\put(2,5){\circle*{0.2}}
\put(3,5){\circle*{0.2}}\put(4,5){\circle*{0.2}}\put(5,5){\circle*{0.2}}\put(6,5){\circle*{0.2}}\put(7,5){\circle*{0.2}}
\put(9,5){\circle*{0.2}}\put(11,5){\circle*{0.2}}
\put(6,7.5){$\downarrow$} \put(-3, 9){\small $x_6=1$}
 \put(-0.2,8.4){\small$1$}\put(0.8,8.4){\small$2$}
\put(2.8,8.4){\small$3$}\put(3.8,8.4){\small$4$}
 \put(5.8,8.4){\small$5$} \put(7.8,8.4){\small$6$}
\put(9.8,8.4){\small$7$}
\qbezier(0,9)(2.5,11.5)(5,9)\qbezier(1,9)(1.5,10.5)(2,9)\qbezier(4,9)(6.5,11.5)(9,9)
\qbezier(6,9)(6.5,10.5)(7,9)
\put(8,9){\circle*{0.2}}\put(10,9){\circle*{0.2}}\put(0,9){\circle*{0.2}}\put(1,9){\circle*{0.2}}\put(2,9){\circle*{0.2}}
\put(3,9){\circle*{0.2}}\put(4,9){\circle*{0.2}}\put(5,9){\circle*{0.2}}\put(6,9){\circle*{0.2}}\put(7,9){\circle*{0.2}}
\put(9,9){\circle*{0.2}}
\put(6,11.5){$\downarrow$} \put(-3, 13){\small $x_5=2$}
 \put(-0.2,12.4){\small$1$}\put(0.8,12.4){\small$2$}
\put(2.8,12.4){\small$3$}\put(3.8,12.4){\small$4$}
 \put(5.8,12.4){\small$5$} \put(7.8,12.4){\small$6$}
\put(8.8,12.4){\small$7$}
\qbezier(0,13)(2.5,15.5)(5,13)\qbezier(1,13)(1.5,14.5)(2,13)
\qbezier(6,13)(6.5,14.5)(7,13)
\put(8,13){\circle*{0.2}}\put(0,13){\circle*{0.2}}\put(1,13){\circle*{0.2}}\put(2,13){\circle*{0.2}}
\put(3,13){\circle*{0.2}}\put(4,13){\circle*{0.2}}\put(5,13){\circle*{0.2}}\put(6,13){\circle*{0.2}}\put(7,13){\circle*{0.2}}\put(9,13){\circle*{0.2}}
\put(4,15.5){$\downarrow$} \put(-3, 17){\small $x_4=0$}
 \put(-0.2,16.4){\small$1$}\put(0.8,16.4){\small$2$}
\put(2.8,16.4){\small$3$}\put(3.8,16.4){\small$4$}
 \put(5.8,16.4){\small$5$} \put(6.8,16.4){\small$6$}
\put(7.8,16.4){\small$7$}
\qbezier(0,17)(2.5,19.5)(5,17)\qbezier(1,17)(1.5,18.5)(2,17)

\put(8,17){\circle*{0.2}}\put(0,17){\circle*{0.2}}\put(1,17){\circle*{0.2}}\put(2,17){\circle*{0.2}}
\put(3,17){\circle*{0.2}}\put(4,17){\circle*{0.2}}\put(5,17){\circle*{0.2}}\put(6,17){\circle*{0.2}}\put(7,17){\circle*{0.2}}

\put(4,19.5){$\downarrow$} \put(-3, 21){\small $x_2=1$}
 \put(-0.2,20.4){\small$1$}\put(0.8,20.4){\small$2$}
\put(2.8,20.4){\small$3$}\put(3.8,20.4){\small$4$}
 \put(4.8,20.4){\small$5$} \put(5.8,20.4){\small$6$}
\put(6.8,20.4){\small$7$}
 \qbezier(1,21)(1.5,22.5)(2,21)

\put(7,21){\circle*{0.2}}\put(0,21){\circle*{0.2}}\put(1,21){\circle*{0.2}}\put(2,21){\circle*{0.2}}
\put(3,21){\circle*{0.2}}\put(4,21){\circle*{0.2}}\put(5,21){\circle*{0.2}}\put(6,21){\circle*{0.2}}

\put(3,23.5){$\downarrow$}
 \put(-0.2,24.4){\small$1$}\put(0.8,24.4){\small$2$}
\put(1.8,24.4){\small$3$}\put(2.8,24.4){\small$4$}
 \put(3.8,24.4){\small$5$} \put(4.8,24.4){\small$6$}
\put(5.8,24.4){\small$7$}

 \put(0,25){\circle*{0.2}}\put(1,25){\circle*{0.2}}\put(2,25){\circle*{0.2}}
\put(3,25){\circle*{0.2}}\put(4,25){\circle*{0.2}}\put(5,25){\circle*{0.2}}\put(6,25){\circle*{0.2}}
\end{picture}
\vspace{-15pt}
\end{center}
\caption{The bijection $\alpha$.} \label{construction1}
\end{figure}
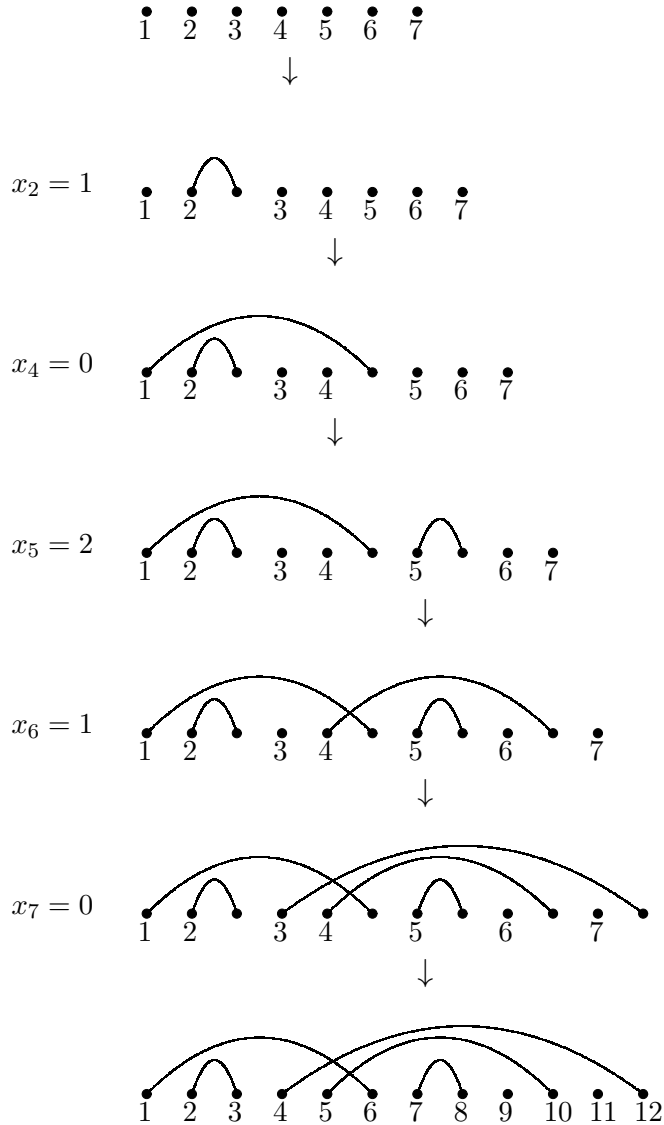

Conversely, given the integer sequence
$x_0x_1x_2x_3x_4x_5x_6x_7=01120210\in \mathcal{S}(8,5)$, the
construction of the corresponding partial matching is illustrated in
Figure \ref{construction1}.

 Combining Theorems   \ref{ths} and \ref{bijection},
 we are led to
 the generating function formula  (\ref{L(n,k)}).

\section{Partitions with no right crossings}
In this section, we present   a bijection between the set
 $\mathcal{P}(n+k-1,k)$ and the set $\mathcal{CT}(n,k)$.  As a result, we
 derive the generating functions for
 partitions avoiding right crossings.  Before we present the bijection, we should recall the notion of $2$-paths defined by Chen
 et al. \cite{chen}. Recall that a pair of two arcs $(i,j)$ and $(j, k)$ with $i<j<k$ in the diagram of a partition  is
 said to be a {\em $2$-path}.

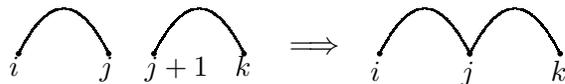
\begin{figure}[h,t]
\begin{center}
\begin{picture}(80,10)
\setlength{\unitlength}{6mm} \linethickness{0.4pt}
 \put(7.8,-0.6){\small$i$}\put(9.8,-0.6){\small$j$} \put(11.8,-0.6){\small$k$}

\qbezier(8,0)(9,2)(10,0)\qbezier(10,0)(11,2)(12,0)
\put(8,0){\circle*{0.1}}\put(10,0){\circle*{0.1}}
\put(10,0){\circle*{0.1}}\put(12,0){\circle*{0.1}}
\put(6,0){$\Longrightarrow$}
\put(-0.2,-0.5){\small$i$}\put(1.8,-0.5){\small$j$}
\put(2.8,-0.5){\small$j+1$} \put(4.8,-0.5){\small$k$}

\qbezier(0,0)(1,2)(2,0)\qbezier(3,0)(4,2)(5,0)
\put(0,0){\circle*{0.1}}\put(2,0){\circle*{0.1}}
\put(3,0){\circle*{0.1}}\put(5,0){\circle*{0.1}}
\end{picture}
\vspace{-15pt}
\end{center}
\caption{change a neighbor alignment to a $2$-path.} \label{fig2}
\end{figure}

\begin{theorem}\label{reduction1}
There is  a bijection between the set $\mathcal{P}(n+k-1, k)$ and
the set $\mathcal{CT}(n,k)$. Moreover, this bijection transforms the
number of  neighbor alignments of a partial matching to the number
of transients of a partition.
\end{theorem}
\pf  Let $M$ be a partial matching in $\mathcal{P}(n+k-1, k)$. we
may reduce it to a partition by the following procedure.
\begin{itemize}
\item  Change   neighbor alignments to
2-paths from left to right until there are no more neighbor
alignments, see Figure \ref{fig2} for an illustration. More
precisely, suppose that there is a neighbor alignment consisting of
two arcs $(i,j)$ and $(j+1, k)$. We  change the arc $(j+1, k)$ to
$(j,k)$ and delete the vertex $j+1$.
\item Delete the
singleton immediately after each closer, except for the last closer.
\item Relabel the vertices   in the natural order.
\end{itemize}
Denote by $P$ the  obtained partition. We claim that there are no
right crossings in the resulting partition $P$. Suppose that there
are two crossing arcs $(i,j)$ and $(k,j+1)$ with $i<k<j$ in the
resulting partition. Since we have either deleted an opener or a
singleton immediately after each closer of $M$, except for the last
closer,   the two crossing arcs forms  a 2-right crossing in the
partial matching $M$, a contradiction. Hence, the claim is proved.

Conversely, given a partition $P'\in \mathcal{CT}(n,k)$, we wish to
construct a partial matching $ M'\in \mathcal{R}(n+k-1, k)$. First,
we add a vertex after each closer, except for the last closer. Then,
we change each pair of arcs $(i,j)$ and $(j,k)$ with $i<j<k$ into a
neighbor alignment. Finally, we relabel the vertices in the natural
order to get the partial matching $M'$.

From the construction of the inverse map, we know that we have
either added a singleton or an opener after each closer of the
partition $P'$. Hence, there are no right crossings and right
nestings in the resulting partial matching. We claim that there are
no 2-right crossings in the resulting partial matching $M'$. Suppose
that there are two crossing arcs $(i,j)$ and $(k,j+2)$ with $i<k<j$
in the resulting partial matching $ M'$. Then $j+1$ is either a
singleton or an opener. From the construction of the reverse map,
the vertex $j+1$ is an added vertex. This implies that the two
crossing arcs is a right crossing in the partition $P'$, a
contradiction. Hence, the claim is proved. This implies that there
is a bijection between the set $\mathcal{P}(n+k-1,k)$ and the set
$\mathcal{CT}(n,k)$. Obviously, this bijection transforms the number
of  neighbor alignments of a partial matching to the number of
transients of a partition. This completes the proof. \qed

 Combining Theorems \ref{L(n,k)} and
\ref{reduction1}, we derive the generating function  Formula
(\ref{T(n,k)}).

\begin{figure}[h,t]
\begin{center}
\begin{picture}(120,10)
\setlength{\unitlength}{5mm} \linethickness{0.4pt}
 \put(-0.2,-0.6){\small$1$}\put(0.8,-0.6){\small$2$}
\put(1.8,-0.6){\small$3$}\put(2.8,-0.6){\small$4$}
 \put(3.8,-0.6){\small$5$} \put(4.8,-0.6){\small$6$}
\put(5.8,-0.6){\small$7$}\put(6.8,-0.6){\small$8$}\put(7.8,-0.6){\small$9$}\put(8.8,-0.6){\small$10$}\put(9.8,-0.6){\small$11$}\put(10.8,-0.6){\small$12$}
\qbezier(0,0)(2.5,2.5)(5,0)\qbezier(1,0)(1.5,1.5)(2,0)\qbezier(3,0)(7,3)(11,0)\qbezier(4,0)(6.5,2.5)(9,0)
\qbezier(6,0)(6.5,1.5)(7,0)
\put(8,0){\circle*{0.1}}\put(10,0){\circle*{0.1}}\put(0,0){\circle*{0.1}}\put(1,0){\circle*{0.1}}\put(2,0){\circle*{0.1}}
\put(3,0){\circle*{0.1}}\put(4,0){\circle*{0.1}}\put(5,0){\circle*{0.1}}\put(6,0){\circle*{0.1}}\put(7,0){\circle*{0.1}}
\put(9,0){\circle*{0.1}}\put(11,0){\circle*{0.1}}
\put(12.5,0){$\Longleftrightarrow$}
\qbezier(15,0)(17,2.5)(19,0)\qbezier(19,0)(19.5,1)(20,0)
\qbezier(16,0)(16.5,1)(17,0)\qbezier(17,0)(19.5,3)(22,0)
\qbezier(18,0)(19.5,2)(21,0) \put(15,0){\circle*{0.1}}
\put(16,0){\circle*{0.1}}\put(17,0){\circle*{0.1}}
\put(18,0){\circle*{0.1}}\put(19,0){\circle*{0.1}}
\put(20,0){\circle*{0.1}}\put(21,0){\circle*{0.1}}\put(22,0){\circle*{0.1}}

 \put(14.8,-0.6){\small$1$} \put(15.8,-0.6){\small$2$} \put(16.8,-0.6){\small$3$}
  \put(17.8,-0.6){\small$4$} \put(18.8,-0.6){\small$5$} \put(19.8,-0.6){\small$6$}
   \put(20.8,-0.6){\small$7$} \put(21.8,-0.6){\small$8$}

\end{picture}
\vspace{-15pt}
\end{center}
\caption{The partial matching $M$ and its corresponding partition
$P$.} \label{fig3}
\end{figure}
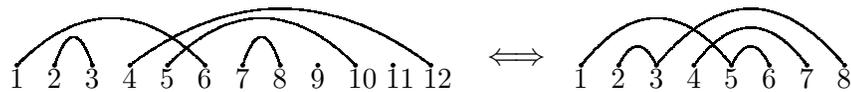

 Figure \ref{fig3} gives an example of a partial matching $M\in
\mathcal{P}(12,5)$ and its corresponding partition $P\in
\mathcal{CT}(8,5)$.

 \vskip 5mm

\noindent{\bf Acknowledgments.}  This work was supported by the
National Natural Science Foundation of China (No. 10901141).


\begin{thebibliography}{100}

\bibitem{melon}
 M. Bousquet-M\'elou, A. Claesson, M. Dukes, S. Kitaev,   $(2+2)$-free posets, ascent sequences and pattern avoiding permutations, {\em
J. Combin. Theory Ser. A} {\bf 117} (2010),  884--909.

\bibitem{chen}
W.Y.C. Chen,  N.J.Y. Fan and A.F.Y.Zhao, Partitions and partial
matchings avoiding neighbor patterns,  {\em Europ.  J. Combin.}, to
appear.

\bibitem{cla}
A. Claesson,  S. Linusson, $n!$ matchings, $n!$ posets, {\em Proc.
Amer. Math. Soc.} {\bf 139} {2011}, 435--449.

\bibitem{Duck}
M. Dukes, R. Parviainen, Ascent sequences and upper triangular
matrices containing non-negative integers, {\em Electron.   J.
combin. } {\bf 17 } (2010),  R53.



\bibitem{levande}
P. Levande, Two new interpretations of the Fishburn numbers and
their refined generating functions, arXiv: math.CO /1006.3013.



 \bibitem{Stoi}
A. Stoimenow, Enumumeration of chord diagrams and an upper bound for
Vassiliev invariants, {\em J. Knot Theory Ramifications}  {\bf 7}
(1998), 93--114.
\bibitem{zag}
D. Zagier, Vassiliev invariants and a stange identity related to the
Dedeking eta-function, {\em Topology}   {\bf 40} (2001), 945--960.


\end{thebibliography}
\end{document}